\documentclass[11pt, a4paper, reqno]{amsart}
\usepackage{amssymb, array, amsmath, amscd, pdfpages, enumerate, amsthm, fixltx2e, setspace, hyperref,mathtools}

\usepackage[utf8]{inputenc}
\usepackage{amssymb}
\usepackage{bm}
\usepackage{graphicx}

\newcommand{\eps}{\varepsilon}

\DeclareMathOperator{\spec}{sp}
\DeclareMathOperator{\spsup}{sp \;sup}

\DeclareMathOperator{\interior}{int}

\newcommand*{\R}{{\mathbb{R}}}

\newcommand*{\N}{{\mathbb{N}}}

\newcommand*{\Abs}[2][default]{\ifthenelse{\equal{#1}{default}}{\left\lvert#2\right\rvert}{\ldelim{#1}{\lvert}#2\rdelim{#1}{\rvert}}}
\newcommand*{\Norm}[2][default]{\ifthenelse{\equal{#1}{default}}{\left\lVert#2\right\rVert}{\ldelim{#1}{\lVert}#2\rdelim{#1}{\rVert}}}

\newcommand*{\Iprod}[3][default]{\ifthenelse{\equal{#1}{default}}{\left\langle#2,#3\right\rangle}{\ldelim{#1}{\langle}#2,#3\rdelim{#1}{\rangle}}}
\newcommand*{\Dualpair}[3][default]{\ifthenelse{\equal{#1}{default}}{\left\langle#2,#3\right\rangle}{\ldelim{#1}{\langle}#2,#3\rdelim{#1}{\rangle}}}

\newcommand*{\ddb}[2][1]{\ifthenelse{\equal{#1}{1}}{\frac{d}{d#2}}{\frac{d^{#1}}{d#2^{#1}}}}
\newcommand*{\pd}[3][1]{\ifthenelse{\equal{#1}{1}}{\frac{\partial{#2}}{\partial{#3}}}{\frac{\partial^{#1}{#2}}{\partial#3^{#1}}}}

\newcommand{\norm}[1]{\|#1\|}

\newcommand*\lenv{{\hbox{\raisebox{-.15ex}{\rotatebox[origin=c]{50}{$\smallsmile$}}\kern-8.65pt\rotatebox[origin=c]{-25}{$\smallsetminus$}}}}
\newcommand*\uenv{{\hbox{\raisebox{-.0ex}{\rotatebox[origin=c]{-45}{$\smallfrown$}}\kern-5.4pt\raisebox{.2ex}{\rotatebox[origin=c]{-5}{\scriptsize\slash}}}}\,\kern+1.5pt}

\newcommand{\sleq}{\preccurlyeq}
\newcommand{\sgeq}{\succcurlyeq}
 
\newcommand*{\Ss}{{\mathcal{S}}}

\newcommand*{\Vv}{{\mathcal{V}}}

\newcommand*{\Aa}{{\mathcal{A}}}
\newcommand*{\Bb}{{\mathcal{B}}}
\newcommand*{\lupp}[1]{\kern0pt^{\kern0pt u \kern0pt}\kern1pt#1}
\newcommand*{\llow}[1]{\kern0pt^{\kern0pt l \kern0pt}\kern1pt#1}
\newcommand*{\rupp}[1]{#1\kern1pt^{\kern0pt u \kern0pt}\kern0pt}
\newcommand*{\rlow}[1]{#1\kern1pt^{\kern0pt l \kern0pt}\kern0pt}
\newcommand*{\ul}[1]{\kern0pt^{\kern0pt u \kern0pt}\kern0pt#1\kern1pt^{\kern0pt l \kern0pt}\kern0pt}
\newcommand*{\lu}[1]{\kern0pt^{\kern0pt l \kern0pt}\kern0pt#1\kern0pt^{\kern0pt u \kern0pt}\kern0pt}
\newcommand*{\bproofname}{Proof}
\newenvironment{bproof}[1][\bproofname]{\begin{proof}[#1]}{\end{proof}}

\newtheorem{thm}{Theorem}[section]
\newtheorem{prop}[thm]{Proposition}
\newtheorem{lemma}[thm]{Lemma}
\newtheorem{cor}[thm]{Corollary}
\theoremstyle{definition}
\newtheorem{defn}[thm]{Definition}

\newtheorem{remark}[thm]{Remark}
\newtheorem{example}[thm]{Example}

\numberwithin{equation}{section}

\begin{document}

\title[Compatible topologies on mixed lattice vector spaces]{Compatible topologies on mixed lattice vector spaces}

\thispagestyle{plain}

\author{Jani Jokela}
\address[J. Jokela]{Mathematics, Faculty of Information Technology and Communication Sciences, 
Tampere University, PO. Box 692, 33101 Tampere, Finland}
\email{jani.jokela@tuni.fi}

\begin{abstract}
A mixed lattice vector space is a partially ordered vector space with two partial orderings, generalizing the notion of a Riesz space.  
Whereas the algebraic theory of mixed lattice structures dates back to the 1970s, 
the topological theory of mixed lattice spaces remains largely unexplored. 
The purpose of this paper is to develop the basic topological theory of mixed lattice spaces.  
A vector topology is said to be compatible with the mixed lattice structure if the mixed lattice operations are continuous. We obtain a characterization of compatible 
mixed lattice topologies, which is similar to the well known Roberts-Namioka characterization of locally solid Riesz spaces.  
Moreover, we study locally convex topologies and the associated seminorms, as well as connections between mixed lattice topologies and locally solid topologies on Riesz spaces. We also briefly discuss asymmetric norms and cone norms on mixed lattice spaces.
\end{abstract}

\subjclass[2010]{%
46A40 
%
}
\keywords{mixed lattice, locally solid, locally full, topological Riesz space, asymmetric norm, cone norm} 

\maketitle

\section{Introduction}
\label{sec:intro}

We recall that a real vector space $\Vv$ together with a partial ordering $\leq$ is called a 
\emph{partially ordered vector space}, if  
\begin{equation*}\label{pospace}
u\leq v \; \implies \; u+w\leq  v+w  \, \quad \textrm{ and } \quad \,
u\leq v \; \implies \; a u\leq a v
\end{equation*}
holds for all $u,v,w\in\Vv$ and $a \in \R_+$. A partially ordered space $\Vv$ is called a \emph{Riesz space}, or \emph{vector lattice}, if $\sup\{x,y\}$ and $\inf \{x,y\}$ exist for all $x,y\in \Vv$. 
 
A mixed lattice vector space is a partially ordered vector space with two partial orderings. 
The precise definition will be given in Section~\ref{sec:sec2}, but the idea is that  
the usual notions of supremum and infimum of two elements in a Riesz space are replaced by asymmetric mixed envelopes which are formed with respect to the two partial orderings. It is then required that these mixed upper and lower envelopes exist for every pair of elements. If the two partial orderings are identical, then the mixed upper and lower envelopes become the usual supremum and infimum, and the mixed lattice vector space is reduced to a Riesz space. In this sense, the concept of a mixed lattice vector space is a generalization of a Riesz space. However, due to the asymmetric behavior of the mixed envelopes, some familiar properties of Riesz spaces (such as the distributive laws and commutativity of the lattice operations) no longer hold in a mixed lattice vector space. In particular, a mixed lattice space is not necessarily a lattice with respect to either of the partial orderings. 

The notion of a mixed lattice semigroup was introduced by Arsove and Leutwiler in 
 connection to their work on axiomatization of potential theory \cite{ars3, ars}.   
The mixed lattice theory in a group setting was later studied by Eriksson--Bique \cite{eri1, eri}. 
More recently, in \cite{jj1} and \cite{jj2} the theory of mixed lattice vector spaces has been developed in a direction that is more parallel with the theory of Riesz spaces, as presented, for example, in \cite{lux}.

All the previous research on mixed lattice structures has focused on the algebraic and order structure, and 
the main purpose of this paper is to develop the topological theory of mixed lattice spaces. The natural starting point is that the topology should be compatible with the mixed lattice structure in the sense that the mixed lattice operations are continuous. This requirement leads to different conditions the topology should satisfy in order to be compatible. A similar topological theory for Riesz spaces and more general ordered vector spaces is well-established (cf. \cite{locsol, cri2, pere}). 

In a general partially ordered vector space the notion of locally full topology provides the most natural setting for topological considerations, and with some modifications, these ideas can be applied also in mixed lattice spaces. We recall here that a subset $S$ of an ordered vector space is called \emph{full} if $x,y\in S$ and $x\leq z\leq y$ imply that $z\in S$. A vector topology is called \emph{locally full} if it has a neighborhood base at zero consisting of full sets.   
On the other hand, in Riesz spaces the fundamental idea is that there is a base of neighborhoods of zero consisting of solid sets. A subset $S$ of a Riesz space is called \emph{solid} if $y\in S$ and $|x|\leq |y|$ imply that $x\in S$. Solid sets are  
defined in terms of the absolute value of an element, which does not necessarily exist in a mixed lattice vector space. 
However, there are ways to generalize the notion of absolute value in mixed lattice spaces, as  
introduced in \cite{jj1} and \cite{jj2}. 
Such generalization can be used to define an analogous notion of a solid set in mixed lattice spaces, resulting in a topological theory that is similar to the theory of locally solid Riesz spaces.

The abovementioned ideas are discussed in Section~\ref{sec:sec3}, where 
we study adaptations of locally full and locally solid topologies in the mixed lattice setting. In Section~\ref{sec:sec35} it is shown that in the locally convex case the topology is determined by a family of seminorms that have additional properties related to the partial orderings.   
Our main results in Sections~\ref{sec:sec3} and \ref{sec:sec35} give fundamental characterizations of mixed lattice topologies, which can be viewed as mixed lattice versions of the well-known Roberts--Namioka theorem in locally solid Riesz spaces. In the general case we obtain a partial result, but in the locally convex case we can give a more complete characterization. 
We also study the connections between mixed lattice topologies and Riesz space topologies in the case that 
the mixed lattice space is a lattice with respect to one of the partial orderings. 
It turns out that a compatible topological structure on a mixed lattice space places some restrictions on the order structure. We show in Section~\ref{sec:sec35} that a finite dimensional normed mixed lattice space is necessarily a lattice with respect to one of the partial orderings. 

In Section~\ref{sec:sec4} we briefly discuss asymmetric norms on mixed lattice spaces. Vector spaces with asymmetric norms are a relatively recent area of research (see \cite{cob} and the references therein), and in mixed lattice spaces asymmetric norms appear naturally due to the asymmetric nature of the mixed envelopes. We close the paper by presenting an application of our results to asymmetric cone norms on $\R^n$. 
More specifically, we show that given any closed cone $C$ in $\R^n$ there is an associated mixed lattice structure on $\R^n$ that gives rise to an asymmetric cone norm corresponding to the cone $C$. 
These type of applications are of interest in other branches of mathematics, 
particularly in convex analysis and related topics \cite{nemeth2020}.

\section{Mixed lattice vector spaces}
\label{sec:sec2}

We begin by recalling some definitions and terminology related to mixed lattice spaces. 
A subset $K$ of a vector space is called a \emph{cone} if \, (i)\;$t K \subseteq K$ for all $t \geq 0$, \, (ii)\;\;$K+K \subseteq K$ and \, (iii) \,$K\cap (-K)=\{0\}$. For any cone $K$ in a vector space there is an associated partial ordering defined by $x\leq y$ iff $y-x\in K$. Then $K$ is called the \emph{positive cone} for the ordering $\leq$.

Suppose next that we have two partial orderings $\leq$ and $\sleq$ on $\Vv$. Here $\leq$ is called the \emph{initial order} and $\sleq$ is called the \emph{specific order}. 
For these two partial orderings $\leq$ and $\sleq$ we define the \emph{mixed upper and lower envelopes} 
\begin{equation}\label{upperenv}
u\uenv v\,=\,\min \,\{\,w\in V: \; w\sgeq u \; \textrm{ and } \; w\geq v \,\}
\end{equation}
and
\begin{equation}\label{lowerenv}
u\lenv v \,=\,\max \,\{\,w\in V: \; w\sleq u \; \textrm{ and } \; w\leq v \,\}, 
\end{equation}
respectively, where the minimum and maximum (whenever they exist) are taken with respect to the initial order $\leq$. These definitions were introduced by Arsove and Leutwiler in \cite{ars3}. 
We observe that these operations are not commutative, i.e. $x\uenv y$ and $y\uenv x$ are not equal, in general. 

A mixed lattice vector space is defined by 
requiring that the mixed envelopes exist for every pair of elements.

\begin{defn}\label{lml}
Let $\Vv$ be a partially ordered real vector space 
with respect to two partial orderings $\leq$ and $\sleq$, and let $\Vv_p$ and $\Vv_{sp}$ be the corresponding positive cones, respectively. Then $(\Vv,\leq,\sleq)$ is called a \emph{mixed lattice vector space}  
if the following conditions hold:
\begin{enumerate}[(1)]
\item
The elements $x\lenv y$ and $x\uenv y$
exist in $\Vv$ for all $x,y\in \Vv$,
\item
$\Vv_{sp}\subseteq \Vv_p$ (i.e. $x\sleq y$ implies $x\leq y$),
\item
$V_{sp}$ is a \emph{mixed lattice cone}, that is, 
the elements $x\uenv y$ and $x\lenv y$ are in $\Vv_{sp}$ whenever $x,y\in \Vv_{sp}$.
\end{enumerate}
\end{defn}

\begin{remark}
A mixed lattice structure can also be defined in a more general setting of an additive group.  
This more general definition 
does not assume the conditions $(2)$ and $(3)$ in the above definition. These conditions are included here because they provide a sufficiently rich structure for developing an interesting theory for mixed lattice spaces. In fact, many important properties of mixed lattice structures depend on these assumptions. For more details on these technicalities, as well as many examples of mixed lattice spaces, we refer to \cite{jj1}. 
\end{remark}

Below we have listed several basic properties of the mixed envelopes that hold in every mixed lattice space. 
For proofs and further discussion on these properties 
we refer to \cite{eri}, \cite{jj1}, and \cite{jj2}. In the following identities and inequalities, 
$x,y,z,u,v$ are elements of a mixed lattice space $\Vv$ and $a\in \R$.

\begin{equation}\label{p0}
x\lenv y \sleq x\sleq x\uenv y \; \textrm{ and } \;  x\lenv y \leq y\leq x\uenv y
\end{equation}
\begin{equation}\label{p1}
x\uenv y \, + \, y \lenv x \, = \, x+y
\end{equation}
\begin{equation}\label{p2}
z \,+ \,x\uenv y \, = \, (x+z)\uenv(y+z) \, \textrm{ and } \, z \,+\, x\lenv y \, = \, (x+z)\lenv(y+z)
\end{equation}
\begin{equation}\label{p3}
x\uenv y \, = \, -(-x \lenv -y)
\end{equation}
\begin{equation}\label{p4}
x\sleq u \; \textrm{ and } \; y\leq v \; \implies \; x\uenv y \leq u\uenv v \; \textrm{ and } \; x
\lenv y \leq  u\lenv v
\end{equation}
\begin{equation}\label{p5a}
x\leq y \, \iff \, y\uenv x = y \, \iff \, x\lenv y = x
\end{equation}
\begin{equation}\label{p5b}
x\sleq y \, \iff \, x\uenv y = y \, \iff \, y\lenv x = x
\end{equation}
\begin{equation}\label{p6}
x\sleq y  \; \implies \; z\uenv x \sleq z\uenv y \; \textrm{ and } \; z\lenv y \sleq  z\lenv y
\end{equation}
\begin{equation}\label{p7}
u\sleq x\sleq z \; \textrm{ and } \; u\sleq y\sleq z \; \implies \; x\uenv y \sleq z \; \textrm{ and } \; u \sleq x\lenv y
\end{equation}
\begin{equation}\label{p8a}
(ax)\lenv(ay) = a(x\lenv y) \; \textrm{ and } \; (ax)\uenv(ay) = a(x\uenv y) \quad  
(a\geq 0)
\end{equation}
\begin{equation}\label{p8b}
(ax)\lenv(ay) = a(x\uenv y) \; \textrm{ and } \;  (ax)\uenv(ay) = a(x\lenv y) \quad  
(a< 0) \smallskip \\
\end{equation}

The following notions of upper and lower parts of an element were introduced in \cite{jj1}. They generalize the concepts of the positive and negative parts of an element in a Riesz space.

\begin{defn}\label{upperparts}
Let $\Vv$ be a mixed lattice vector space and $x\in \Vv$. The elements $^{u}x=x\uenv 0$ and $^{l}x=(-x)\uenv 0$ are called the \emph{upper part} and \emph{lower part} of $x$, respectively. Similarly, the elements $x^{u}=0\uenv x$ and $x^{l}=0\uenv (-x)$ are called \emph{specific upper part} and \emph{specific lower part} of $x$, respectively. The elements \, $\ul{x}= x\uenv (-x)$ \, and \, $\lu{x}=(-x)\uenv x$ \, are called the \emph{(asymmetric) generalized absolute values} of $x$.
\end{defn}

From the above definitions we observe that for the specific upper and lower parts we have $x^{u}\sgeq 0$ and $x^{l}\sgeq 0$, and for the upper and lower parts $^{u}x\geq 0$ and $^{l}x\geq 0$.

The following useful result is from \cite[Theorem 2.12]{jj1}.

\begin{thm}[{\cite[Theorem 2.12]{jj1}}]\label{mlg_ehto}
Let $(\Vv,\leq,\sleq)$ be a partially ordered vector space with two partial orderings. Then $\Vv$ is a mixed lattice vector space if and only if one of the elements \,$\lupp{x}$,\, $\rupp{x}$,\, $\llow{x}$ or $\rlow{x}$\, exists for all $x\in \Vv$.
\end{thm}

The upper and lower parts and the generalized absolute values 
have several important basic properties, which were proved in \cite{jj1}. These properties are given in the next theorem.

\begin{thm}\label{absval}
Let $\Vv$ be a mixed lattice vector space and $x\in \Vv$. Then the following hold.
\begin{enumerate}[(a)]
\item
\; $\lupp{x}\,=\,\llow{(-x)}$ \quad \textrm{and} \quad $\rupp{x}\,=\,\rlow{(-x)}$. 
\item
\; $x \, = \, \rupp{x} \, -\,\llow{x} \, = \, \lupp{x} \, - \, \rlow{x}$.
\item
\;  $\ul{x} \, =\, \lupp{x}\uenv \rlow{x} \, =\, \lupp{x} \,+\,\rlow{x}$ \quad \textrm{and} \quad $\lu{x} \, =\,\llow{x}\uenv \rupp{x} \, =\, \llow{x} \,+\,\rupp{x}$.
\item
\; $\ul{x} \,=\,\lu{(-x)}$.
\item
\;
$\lupp{x}+\lupp{y} \geq \lupp{(x+y)}$, \quad $\rlow{x}+\rlow{y} \geq \rlow{(x+y)}$ \quad and \quad $\ul{x}+\ul{y} \geq \ul{(x+y)}$.
\item
\;
$\rupp{x}+\rupp{y} \geq \rupp{(x+y)}$, \quad $\llow{x}+\llow{y} \geq \llow{(x+y)}$ \quad and \quad $\lu{x}+\lu{y} \geq \lu{(x+y)}$.
\item
\; $\rupp{x}\lenv \llow{x} \, = \, 0\, = \, \rlow{x}\lenv \lupp{x}$.
\item
\; $\rupp{x}\uenv \llow{x}\, = \,\lupp{x}\,+\,\llow{x}\, = \,\rlow{x}\,+\,\rupp{x}\, = \,\rlow{x}\uenv \lupp{x}$
\item
\; $x\sgeq 0$ \, if and only if \, $x=\lu{x}=\ul{x}=\lupp{x}=\rupp{x}$\, and \,$\llow{x}=\rlow{x}=0$.
\item
\; $x\geq 0$ \, if and only if \, $x=\ul{x}=\lupp{x}$\, and \,$\rlow{x}=0$.
\item
\; $\ul{x}\geq 0$\,  and  \,$\lu{x}\geq 0$. Moreover,  $\ul{x}=\lu{x}= 0$\,  if and only if  $x=0$.
\item
\; $\ul{(ax)}=a\,\ul{x}$  \, and \, $\lu{(ax)}=a\,\lu{x}$ \, \textrm{for all} $a\geq 0$.
\item
\; $\ul{(ax)}=|a|\,\lu{x}$  \, and \, $\lu{(ax)}=|a|\,\ul{x}$ \, \textrm{for all} $a<0$.
\end{enumerate}
\end{thm}

We should note here that in a mixed lattice vector space $\Vv$ the positive cone $\Vv_p=\{x\in \Vv: x\geq 0\}$ is always generating. This is a consequence of Theorem \ref{absval}(b).

It is possible to combine the two asymmetric generalized absolute values to define a symmetric version of the generalized absolute value that retains most of the important properties of the absolute value. This notion was introduced in \cite{jj2}, and it will be essential for describing topologies on mixed lattice vector spaces.

\begin{defn}\label{symmap}
Let $\Vv$ be a mixed lattice vector space and $x\in \Vv$. The element \, $s(x)=\frac{1}{2}(\ul{x}+\lu{x})$\, is called the \emph{symmetric generalized absolute value} of $x$. 
\end{defn}

The next theorem lists the most important properties of the symmetric absolute value. 
They are rather straightforward consequences of Theorem \ref{absval} (for proofs, see \cite[Theorem 5.17]{jj2}).

\begin{thm}[{\cite[Theorem 5.17]{jj2}}]\label{sav}  
Let $\Vv$ be a mixed lattice vector space and $x\in \Vv$. Then the following hold.
\begin{enumerate}[(a)]
\item
\; $s(x)\,=\,\rupp{x}\uenv \llow{x}\, = \,\lupp{x}\,+\,\llow{x}\, = \,\rlow{x}\,+\,\rupp{x}\, = \,\rlow{x}\uenv \lupp{x}$.
\item
\; $s(\alpha x)=|\alpha|s(x)$ for all $\alpha\in\R$. 
\item
\; $s(x)\sgeq 0$ \, and \, $s(x)\geq 0$. Moreover,  $s(x)= 0$  if and only if  $x=0$.
\item
\; $x\sgeq 0$ \, if and only if \, $x=s(x)=\rupp{x}$. 
\item
\; $s(s(x))=s(x)$.
\item
\;
$s(x+y)\leq s(x)+s(y)$ 
\end{enumerate}
\end{thm}

The following result will also be useful in the next section.

\begin{lemma}\label{sav_lemma}
If $\Vv$ is a  mixed lattice vector space then 
$s(\rupp{x}-\rupp{y})\leq s(x-y)$ for all $x,y\in \Vv$.
\end{lemma}

\begin{bproof}
Since $x=y+(x-y)$, Theorem \ref{absval}(f) implies that $\rupp{x} \leq \rupp{y}+\rupp{(x-y)}$ and thus $\rupp{x}- \rupp{y}\leq \rupp{(x-y)}$. Since $0\sleq 0$, we have from \eqref{p4} and Theorem \ref{absval}(i) that
$$
\rupp{(\rupp{x}- \rupp{y})} = 0 \uenv (\rupp{x}- \rupp{y}) \leq 0 \uenv \rupp{(x-y)} = \rupp{(\rupp{(x-y)})} =\rupp{(x-y)}.
$$
Exchanging $x$ and $y$ gives similarly 
$\rupp{(\rupp{y}- \rupp{x})}\leq \rupp{(y-x)}$. But by Theorem \ref{absval}(a) we have $\rupp{(\rupp{y}- \rupp{x})}=\rlow{(\rupp{x}- \rupp{y})}$ \, and \, $\rupp{(y-x)}=\rlow{(x-y)}$, and hence $\rlow{(\rupp{x}- \rupp{y})}\leq \rlow{(x-y)}$. Adding the two inequalities gives
$$
\rupp{(\rupp{x}- \rupp{y})} + \rlow{(\rupp{x}- \rupp{y})} \leq \rupp{(x-y)} + \rlow{(x-y)},
$$ 
or equivalently, \, $s(\rupp{x}- \rupp{y})\leq s(x-y)$, by Theorem \ref{sav}(a).
\end{bproof}


\section{Topological mixed lattice vector spaces}
\label{sec:sec3}


We collect some basic facts and terminology from the theory of topological vector spaces. A more detailed account is given in \cite{cri2} and \cite{tvs2}. 
Let $S$ be a subset of a vector space $\Vv$. Then $S$ is called \emph{symmetric}, if $S=-S$, and $S$ is called \emph{balanced} if $x\in S$ implies $ax\in S$ whenever $|a|\leq 1$. Moreover, 
$S$ is called \emph{absorbing} if for every $x\in \Vv$ there exists a number $r>0$ such that $x\in aS$ whenever $|a|\geq r$, and $S$ is called \emph{convex} if $ax+(1-a)y\in S$ for all $x,y\in S$ and $a\in [0,1]$.

If $\tau$ is a topology on $\Vv$ such that vector addition and scalar multiplication are continuous mappings then $\tau$ is called a \emph{vector topology} and $(\Vv, \tau)$ is a \emph{topological vector space}. 
A collection $\Bb$ of neighborhoods of zero of $\Vv$ is called a \emph{neighborhood base at zero} if for every neighborhood $U$ of zero there exists $V\in \Bb$ such that $V\subseteq U$. 
Every vector topology is translation invariant, that is, for any $y\in \Vv$ the map $x\mapsto x+y$ is a homeomorphism.
This implies that if $x\in \Vv$ and $\Bb$ is a neighborhood base at zero, then $x+\Bb=\{x+B:b\in\Bb\}$ is a neighborhood base at $x$.
This means that the topology is completely determined by the neighborhoods of zero, and for this reason we usually only need to work with neighborhoods of zero.

Every topological vector space has a base at zero consisting of balanced absorbing sets. Furthermore, if $U$ is any neighborhood of zero, then there is a balanced neighborhood $V$ of zero such that $V+V\subseteq U$. A vector topology $\tau$ is called \emph{locally convex}, if $\tau$ has a neighborhood base at zero consisting of convex sets. In this case, for any neighborhood $U$ of zero there exists a convex, balanced and absorbing neighborhood $V$ of zero such that $V\subseteq U$. 

A mapping $p:\Vv \to \R$ is called a \emph{seminorm} if 
$(i)$ $p(x)\geq 0$ for all $x\in \Vv$, \, 
$(ii)$ $p(ax)= |a|p(x)$ for all $x\in \Vv$, $a\in \R$, and 
$(iii)$ $p(x+y)\leq p(x)+p(y)$ for all $x,y\in \Vv$.   
A locally convex vector topology is determined by a family of seminorms. We may always assume that the generating family $\mathcal{P}$ of seminorms is directed, i.e. for every $p_1,p_2\in \mathcal{P}$ there exists $p\in \mathcal{P}$ such that $p(x)\geq \max \{p_1(x),p_2(x)\}$ for all $x\in \Vv$. In this case, a base of neighborhoods of zero is given by the sets $V(p,\eps)=\{x:p(x)< \eps\}$ where $p\in\mathcal{P}$ and $\eps>0$. 
For an absorbing subset $U$ the mapping defined as $p_{\scriptscriptstyle U}(x)=\inf \{t>0:x\in tU\}$ for each $x\in \Vv$ is called the \emph{Minkowski functional} of $U$. If $U$ is an absorbing, convex and balanced set then $p_{\scriptscriptstyle U}$ is a seminorm. 

Let $\Vv$ be a topological vector space. A map $f:\Vv\to \Vv$ is called \emph{uniformly continuous}, if for every neighborhood $V$ of zero 
there exists a neighborhood $U$ of zero 
such that $x-y\in U$ implies $f(x)-f(y)\in V$. For example, the inversion map $x\mapsto -x$, addition and translation are all uniformly continuous. Uniform continuity implies continuity, and the composition of uniformly continuous maps is uniformly continuous.

If $(\Vv,\leq)$ is an ordered vector space then a subset $U$ is called \emph{full} if $x,y\in U$ and $x\leq z\leq y$ imply that $z\in U$. 
An ordered vector space equipped with a vector topology $\tau$ is called an \emph{ordered topological vector space} if $\tau$ is a \emph{locally full} topology, i.e. it has a neighborhood base at zero consisting of full sets. A vector topology is locally full if and only if every neighborhood $U$ of zero contains a neighborhood $V$ of zero such that $y\in V$ and $0\leq x\leq y$ imply that $x\in V$. 
For more information on ordered topological vector spaces, see \cite{cri2} and \cite{pere}.

We will now turn our attention to topologies on a mixed lattice space. 
In order to obtain useful results the topology should be compatible with the mixed lattice structure, and this motivates the following definition.

\begin{defn}
Let $\Vv$ be a mixed lattice vector space equipped with a vector space topology $\tau$. Then $\tau$ is called a \emph{mixed lattice topology} and $(\Vv,\tau)$ is called a \emph{topological mixed lattice space} if the mixed upper and lower envelopes are continuous mappings from $\Vv\times \Vv$ to $\Vv$ (here it is understood that $\Vv\times \Vv$ is equipped with the product topology).
\end{defn}

In Riesz spaces the compatibility is achieved by requiring that the topology is \emph{locally solid}, i.e. it has a local base consisting of solid sets. In fact, the following well-known Roberts--Namioka theorem holds in topological Riesz spaces. 
The theorem is due to G. T. Roberts, who introduced the notion of a locally solid topology in \cite{rob}, and I. Namioka, who later extended the theorem in \cite{nam}.

\begin{thm}[Roberts--Namioka]\label{rob_nam}
Let $(\Vv,\tau)$ be a 
Riesz space equipped with a vector topology $\tau$. 
The following conditions are equivalent.
\begin{enumerate}[(a)]
\item
$\tau$ is a locally solid topology.
\item
$\tau$ is locally full and the lattice operations are continuous at zero.
\item
The lattice operations are uniformly continuous.
\end{enumerate}
\end{thm}

As one of the the main results of this paper we will show that a similar characterization holds in topological mixed lattice spaces. 
Let us begin by considering the uniform continuity of the mixed lattice operations.

\begin{prop}\label{unicont}
Let $(\Vv,\tau)$ be a  
mixed lattice vector space equipped with a vector topology $\tau$. Then the following conditions are mutually equivalent.
\begin{enumerate}[(a)]
\item
The map $(x,y)\mapsto x\uenv y$ from $\Vv\times \Vv$ to $\Vv$ is uniformly continuous.
\item
The map $(x,y)\mapsto x\lenv y$ from $\Vv\times \Vv$ to $\Vv$ is uniformly continuous.
\item
The map $x\mapsto \lupp{x}$ from $\Vv$ to $\Vv$ is uniformly continuous.
\item
The map $x\mapsto \llow{x}$ from $\Vv$ to $\Vv$ is uniformly continuous.
\item
The map $x\mapsto \lu{x}$ from $\Vv$ to $\Vv$ is uniformly continuous.
\item
The map $x\mapsto \ul{x}$ from $\Vv$ to $\Vv$ is uniformly continuous.
\item
The map $x\mapsto \rlow{x}$ from $\Vv$ to $\Vv$ is uniformly continuous.
\item
The map $x\mapsto \rupp{x}$ from $\Vv$ to $\Vv$ is uniformly continuous.
\end{enumerate}
\end{prop}

\begin{bproof}
\emph{(a)}$\implies$\emph{(b)} \, This follows from the identity $x\lenv y  =  -(-x \uenv -y)$ and uniform continuity of the map $x\mapsto -x$. \\
\emph{(b)}$\implies$\emph{(c)} \, Follows from the identity \, $\lupp{x}=-((-x)\lenv 0)$.\\
\emph{(c)}$\implies$\emph{(d)} \, Follows from the identity \, $\llow{x}=\lupp{(-x)}$. \\
\emph{(d)}$\implies$\emph{(e)} \, This follows since \, $\lu{x}=2\llow{x}+x$ and addition and scalar multiplication are uniformly continuous. \\
\emph{(e)}$\implies$\emph{(f)} \, Follows from the identity \, $\ul{x}=\lu{(-x)}$. \\
\emph{(f)}$\implies$\emph{(g)} \, Follows from the identity \, $2\rlow{x}=\ul{x}-x$. \\ 
\emph{(g)}$\implies$\emph{(h)} \, Follows from the identity \, $\rupp{x}=\rlow{(-x)}$. \\
\emph{(h)}$\implies$\emph{(a)} \, If $x\mapsto \rupp{x}$ is uniformly continuous then, since addition is  uniformly continuous, it follows that 
$$
x\uenv y = x + 0\uenv(y-x) = x+ \rupp{(y-x)}
$$  
is also uniformly continuous as a composition of uniformly continuous mappings.
\end{bproof}

\begin{remark}\label{huom1}
In any topological mixed lattice space the mapping $x\mapsto s(x)$ is also continuous. This follows from the identities in Theorem \ref{sav}(a). 
\end{remark}

We will next define different types of sets that will be used in the study of mixed lattice topologies. 

\begin{defn}
Let $A$ be a subset of $\Vv$ and define the sets 
$$
MF_{1}(A)=\{y\in \Vv:x\sleq y \leq z \textrm{  for some  } x,z\in A\}
$$ 
and 
$$
MF_{2}(A)=\{y\in \Vv:x\leq y \sleq z \textrm{  for some  } x,z\in A\}.
$$ 
The set $MF_{1}(A)$ is called the \emph{type 1 mixed-full hull} of $A$, and $MF_{2}(A)$ is called the \emph{type 2 mixed-full hull} of $A$. 
If $A=MF_{1}(A)$ then $A$ is called a \emph{type 1 mixed-full set}, and if $A=MF_{2}(A)$ then $A$ is called a \emph{type 2 mixed-full set}. 
\end{defn}

The mixed-full hulls of a subset $A$ are the smallest mixed-full sets containing $A$. Alternatively, these sets are given by $MF_{1}(A)=(A+\Vv_{sp})\cap (A-\Vv_p)$ and $MF_{2}(A)=(A+\Vv_{p})\cap (A-\Vv_{sp})$. 

Next we show that with regards to neighborhood bases at zero the two different notions of mixed-full  
sets are equivalent, and so for all practical purposes, we only need to consider one type of mixed-full 
neighborhoods.

\begin{prop}\label{mixed-full-char}
Let $(\Vv,\tau)$ be a 
mixed lattice vector space equipped with a vector topology $\tau$. 
The following conditions are equivalent.
\begin{enumerate}[(a)]
\item
Each neighborhood of zero contains a type 2 mixed-full neighborhood of zero.
\item
Each neighborhood of zero contains a neighborhood $V$ of zero with the following property: if $y\in V$ and $0\leq x \sleq y$ then $x\in V$.
\item
Each neighborhood of zero contains a type 1 mixed-full neighborhood of zero.
\item
Each neighborhood of zero contains a neighborhood $V$ of zero with the following property: if $y\in V$ and $0\sleq x \leq y$ then $x\in V$.
\end{enumerate}
\end{prop}

\begin{bproof}
The implication $(a)\implies (b)$ is obvious. 
Suppose that $(b)$ holds. Let $U$ be a neighborhood of zero and let $W_1$ be a neighborhood of zero such that $W_1+W_1\subseteq U$. Let $W_2\subseteq W_1$ be a 
neighborhood of zero with the property given in $(b)$. Choose a balanced neighborhood $V$ of zero such that $V+V\subseteq W_2$. If $x\in MF_1(V)$ then there exist $y,z\in V$ such that $y\sleq x\leq z$. Then $0\leq z-x\sleq z-y$ where $z-y\in W_2$, and it follows that $z-x\in W_2\subseteq W_1$. Consequently, since $z\in V\subseteq W_1$, we have $x=z-(z-x)\in W_1+W_1\subseteq U$. Hence, $MF_1(V)\subseteq U$ and so $(c)$ holds. The implication $(c)\implies (d)$ is obvious and the implication $(d)\implies (a)$ is proved in a similar way as $(b)\implies (c)$, and so 
the proof is complete.
\end{bproof}

The above results motivate the following definition.

\begin{defn}
Let $\Vv$ be a mixed lattice vector space with a vector topology $\tau$. 
Then $\tau$ is called \emph{locally mixed-full} if $\tau$ satisfies the equivalent conditions of Proposition \ref{mixed-full-char}.  
\end{defn}

The equivalence of the conditions given in Proposition \ref{mixed-full-char} also motivate the following convention for simplifying the terminology.  
We shall say that 
a neighborhood $V$ of zero is \emph{mixed-full} if $x\in V$ whenever $0\sleq x \leq y$ or $0\leq x \sleq y$ holds with $y\in V$.

The following basic result provides more information about mixed-full sets and locally mixed-full topologies. 

\begin{prop}\label{nbhoods} 
Let $(\Vv,\tau)$ be a topological mixed lattice space.
\begin{enumerate}[(a)]
\item
If $A$ is a balanced and absorbing set then $MF_1(A)=-MF_2(A)$ and the set $B=MF_1(A)+ MF_2(A)$ is balanced, absorbing and mixed-full (i.e. it has the properties $(b)$ and $(d)$ of Proposition \ref{mixed-full-char}). 
If, in addition, $A$ is convex then $MF_1(A)$, $MF_2(A)$ and $B$ are convex too. 
\item 
If $\tau$ is a locally mixed-full topology then every neighborhood of zero contains a balanced absorbing neighborhood of zero with the properties $(b)$ and $(d)$ of Proposition \ref{mixed-full-char}. If $\tau$ is also locally convex then every neighborhood of zero contains a convex, balanced and absorbing neighborhood of zero with the aforementioned properties.
\end{enumerate}
\end{prop}

\begin{bproof}
$(a)$\,  First we note that if $A$ is absorbing, then so are $MF_1(A)$ and $MF_2(A)$ since $A$ is contained in these sets. Let $A$ be a balanced set. If $x\in MF_1(A)$ then $y\sleq x\leq z$ for some $y,z\in A$. This implies that $-z\leq -x\sleq -y$ where $-y,-z\in A$ since $A$ is symmetric. Hence $-x\in MF_2(A)$. The reverse inclusion is similar, so $MF_1(A)=-MF_2(A)$. Consequently, $B=MF_1(A)+ MF_2(A)$ is a symmetric set. 
Let $x\in B$ and $0\leq t\leq 1$. Then $x=x_1+x_2$ where $x_1\in MF_1(A)$ and $x_2\in MF_2(A)$. If $x_1\in MF_1(A)$ then $y\sleq x_1\leq z$ for some $y,z\in A$. This implies that  $ty\sleq tx_1\leq tz$, and since $A$ is balanced, we have $ty,tz\in A$, and the last inequality shows that $tx_1\in MF_1(A)$. A similar argument shows that $tx_2\in MF_2(A)$, and so $tx\in B$. This shows that $B$ is balanced.    
In addition, if $A$ is convex then $MF_{1}(A)=(A+\Vv_{sp})\cap (A-\Vv_p)$ and $MF_{2}(A)=(A+\Vv_{p})\cap (A-\Vv_{sp})$ are also convex, as the sum and intersection of convex sets are again convex. It follows that $B$ is also convex. 
It remains to show that $B$ has the stated properties. For this, let $0\sleq x\leq y$ where $y\in B$. Then $y=y_1+y_2$ where $y_1\in MF_1(A)$ and $y_2\in MF_2(A)$. From the above inequality we get $-y_2\sleq x-y_2\leq y_1$, where $y_1\in MF_1(A)$ and $-y_2\in -MF_2(A)=MF_1(A)$. Hence, $x-y_2\in MF_1(A)$ and so $x=(x-y_2)+y_2 \in MF_1(A)+MF_2(A)=B$. A similar argument shows that $B$ has the property $(d)$ of Proposition \ref{mixed-full-char}.

$(b)$\, Let $\tau$ be a locally mixed-full topology. Let $U$ be a neighborhood of zero and choose a balanced neighborhood $W_0$ of zero and a type 1 mixed-full neighborhood $W_1$ of zero such that $W_0+W_0\subseteq U$ and $W_1\subseteq W_0$. Next, let $V$ be a balanced absorbing neighborhood of zero such that $V\subseteq W_1$. Then $MS_1(V)\subseteq W_1\subseteq W_0$, and since $W_0$ is balanced, we also have $MF_2(V)=-MF_1(V)\subseteq W_0$. It follows that $B=MF_1(V)+ MF_2(V)\subseteq W_0+W_0\subseteq U$, where $B$ is balanced and absorbing, and it has the stated properties, by part $(a)$. If $\tau$ is also locally convex then the set $V$  can be chosen to be convex and balanced. The desired result then follows by part $(a)$.
\end{bproof}

For the characterization of mixed lattice topologies we need a few more definitions.

\begin{defn}
Let $(\Vv)$ be a mixed lattice vector space. 
\begin{enumerate}[(i)]
\item
A subset $A\subseteq \Vv$ is called \emph{symmetric-solid} if $x\in A$ and $s(y)\leq s(x)$ together imply that $y\in A$. 
The set 
$SH(A)=\{y\in \Vv: s(y) \leq s(x) \textrm{  for some  } x\in A\}$ 
is called the \emph{symmetric-solid hull} of $A$. 
\item
A subset $S\subseteq \Vv$ is called \emph{($\leq$)-full} if $x,y\in S$ and $y\leq z\leq x$ together imply that $z\in S$. Similarly, a subset $S\subseteq \Vv$ is called \emph{($\sleq$)-full} if $x,y\in S$ and $y\sleq z\sleq x$ together imply that $z\in S$.
\end{enumerate}
\end{defn}

\begin{remark}\label{ssolid_remark}
We observe that $SH(A)$ is the smallest symmetric-solid set containing $A$ (with respect to set inclusion).  
Moreover, it follows easily from the properties given in Theorem \ref{sav} that every symmetric-solid set is symmetric, and if $A$ is balanced and absorbing, then so is $SH(A)$. 
\end{remark}

\begin{defn} 
A vector topology $\tau$ on a mixed lattice vector space is called \emph{locally symmetric-solid} if $\tau$ has a base at zero consisting of symmetric-solid sets. 
A vector topology $\tau$ is called \emph{locally ($\leq$)-full} if $\tau$ has a base at zero consisting of ($\leq$)-full sets. Similarly, $\tau$ is called \emph{locally ($\sleq$)-full} if $\tau$ has a base at zero consisting of ($\sleq$)-full sets. 
\end{defn}

The next theorem is the first version of our main result on mixed lattice topologies.

\begin{thm}\label{locsymmsol2}
Let $(\Vv,\tau)$ be a 
mixed lattice vector space equipped with a vector topology $\tau$. 
Consider the following statements:
\begin{enumerate}[(a)]
\item
$\tau$ is locally symmetric-solid.
\item
$\tau$ is locally ($\leq$)-full and the mixed lattice operations are continuous at zero.
\item
The mixed lattice operations are uniformly continuous.
\item
$\tau$ is locally mixed-full and the mixed lattice operations are continuous at zero.
\end{enumerate}
Then the statements (a) and (b) are equivalent, and they both imply (c), and (c) implies (d).
\end{thm}

\begin{bproof}
\emph{(a)}$\implies$\emph{(b)} \,
We will first show that the condition \emph{(a)} implies uniform continuity of the map $x\mapsto \rupp{x}$. 
Let $U$ be a neighborhood of zero. Then there exists a symmetric-solid neighborhood $V$ of zero such that $V\subseteq U$. If $x-y\in V$, then $s(x-y)\in V$ and by Lemma \ref{sav_lemma} we have
$s(\rupp{x}-\rupp{y})\leq s(x-y)$
and since $V$ is symmetric-solid we have $\rupp{x}-\rupp{y}\in V\subseteq U$. This shows that the map $x\mapsto \rupp{x}$ is uniformly continuous, and so by Proposition \ref{unicont} the mixed lattice operations are uniformly continuous. In particular, they are continuous at zero. 

To show that $\tau$ is locally $(\leq)$-full we only need to show that every symmetric-solid set is $(\leq)$-full. For this, let $V$ be a symmetric-solid set with $x\in V$ and $0\leq y\leq x$. Then $\rlow{y}=\rlow{x}=0$ and $\rupp{y}=0\uenv y\leq 0\uenv x=\rupp{x}$, by \eqref{p4}. Thus, 
$$
s(y)=\rupp{y}+\rlow{y}=\rupp{y}\leq\rupp{x}=\rupp{x}+\rlow{x}=s(x).
$$
Since $V$ is symmetric-solid, we deduce that $y\in V$. Hence, $V$ is $(\leq)$-full.

\emph{(b)}$\implies$\emph{(a)} \,
Let $U$ be a neighborhood of zero and choose a $(\leq)$-full and balanced neighborhood $V$ of zero such that $V\subseteq U$. Choose another neighborhood $W$ of zero such that $W+W\subseteq V$. Since the map $x\mapsto\lupp{x}$ is continuous we can choose a balanced neighborhood $V_1$ of zero such that $x\in V_1$ implies $\lupp{x}\in W$. Since $V_1$ is balanced we have $-x\in V_1$ and so $\lupp{(-x)}=\llow{x}\in W$. Then $s(x)=\lupp{x}+\llow{x}\in W+W\subseteq V$. Now, if $s(y)\leq s(x)$ then we have the following inequalities
$$
-s(x)\leq -s(y)\leq -\rlow{y}\leq y\leq\lupp{y}\leq s(y)\leq s(x).
$$
Since $V$ is balanced we have $\pm s(x)\in V$, and so $y\in V$ since $V$ is $(\leq)$-full. Hence $U$ contains the symmetric-solid hull of $V_1$ and so \emph{(a)} holds. 

The implication \emph{(a)}$\implies$\emph{(c)} was already proved. 
To finish the proof, we will show that uniform continuity of the mixed lattice operations implies the statement $(d)$.  If the mixed lattice operations are uniformly continuous, then they are certainly  continuous at zero. 
Let $U$ be a neighborhood of zero and choose a balanced neighborhood  
$W$ of zero such that $W+W\subseteq U$. Since the maps $x\mapsto\llow{x}$ and $x\mapsto\rupp{x}$ are uniformly continuous we can choose a 
neighborhood $V$ of zero such that $z-y\in V$ implies $\rupp{z}-\rupp{y}\in W$ and $\llow{z}-\llow{y}\in W$. 
Assume that $u\sleq x \leq v$ with $u,v\in V$. We write $v=x-(x-v)$ and $u=x-(x-u)$ and note that $x-v\leq 0$ and $x-u\sgeq 0$ imply that $\rupp{(x-v)}=0$ and $\llow{(x-u)}=0$. It follows that $\rupp{x}=\rupp{x}-\rupp{(x-v)}\in W$ and $\llow{x}=\llow{x}-\llow{(x-u)}\in W$. 
Hence $x=\rupp{x}-\llow{x}\in W + W\subseteq U$. This shows that $MF_1(V)$  
is contained in $U$. Hence, the condition $(d)$ holds and the proof is complete.
\end{bproof}

Since every mixed-full set is $(\sleq)$-full, 
the preceding theorem show in particular that locally symmetric-solid mixed lattice spaces are always ordered topological vector spaces with respect to both partial orderings, and therefore all the well-known results for ordered topological vector spaces 
hold in mixed lattice spaces. For more on 
these results, see 
 \cite{cri2} and \cite{pere}. 

We do not know if all the conditions in Theorem \ref{locsymmsol2} are equivalent in general. However, under certain stronger assumptions we obtain a more complete characterization, which will be given later (c.f. Theorem \ref{complete}). First we need to introduce some additional concepts.


Let $A$ be a subset of $\Vv$ and define the sets 
$$
MS_{1}(A)=\{y\in \Vv: -s(x)\sleq y \leq s(x) \textrm{  for some  } x\in A\}
$$ 
and 
$$
MS_{2}(A)=\{y\in \Vv: -s(x)\leq y \sleq s(x) \textrm{  for some  } x\in A\}.
$$

These sets are useful technical devices in the study of topological properties of mixed lattice spaces. 
We should point out that 
the set $A$ is not usually contained in $MS_{1}(A)$ 
or $MS_{2}(A)$. 
However, it is clear that if $A$ is any non-empty set then $0\in MS_1(A)$ and $0\in MS_2(A)$. 
Other basic properties of these sets are given in the next proposition.

\begin{prop}\label{msolid_prop}
Let $\Vv$ be a mixed lattice space. 
\begin{enumerate}[(a)]
\item
If $x\in A$ then $\pm s(x)\in MS_1(A)$ and $\pm s(x)\in MS_2(A)$. Moreover, $MS_{1}(A)=-MS_{2}(A)$.
\item
If the cone $\Vv_{sp}$ is generating and $A$ is an absorbing set then $MS_{1}(A)$ and $MS_{2}(A)$ are also absorbing. 
\item
If $0\leq t \leq 1$ and $x\in MS_1(A)$ then $tx\in MS_1(A)$, and similarly, if $x\in MS_2(A)$ then $tx\in MS_2(A)$. Consequently, the set $MS_1(A)\cup MS_2(A)$ is balanced.
\end{enumerate}
\end{prop}

\begin{bproof}
$(a)$\, 
If $x\in A$ then, in particular, $-s(x)\sleq \pm s(x)\leq s(x)$ holds. 
This implies that $\pm s(x)\in MS_1(A)$. Similarly, $\pm s(x)\in MS_2(A)$. 
If $x\in MS_{1}(A)$ then $-s(y)\sleq x\leq s(y)$ holds for some $y\in A$. This is equivalent to $-s(y)\leq -x\sleq s(y)$ for some $y\in A$, and so $-x\in MS_{2}(A)$. 
Hence $MS_{1}(A)=-MS_{2}(A)$.

$(b)$\,  
Let $A$ be an absorbing set. If $x\in \Vv$ and the cone $\Vv_{sp}$ is generating then $x=u-v$ for some $u,v\in \Vv_{sp}$. Since $A$ is absorbing, there exists some $t>0$ such that $t(u+v)\in A$, and since $t(u+v)\sgeq 0$, we have $s(t(u+v))=t(u+v)$ and $-t(u+v)\sleq tx\sleq t(u+v)$. Thus $-s(t(u+v))\sleq tx \sleq s(t(u+v))$, and this implies that $tx\in MS_1(A)$ and $tx\in MS_2(A)$. Hence, $MS_1(A)$ and $MS_2(A)$ are absorbing.

$(c)$\, 
Assume first that $0\leq t\leq 1$ and $x\in MS_1(A)$ for some set $A$. Then there exist some $y\in A$ such that $-s(y)\sleq x\leq s(y)$. Since $s(y)\sgeq 0$, this implies that $-s(y)\sleq -ts(y)\sleq tx\leq ts(y)\leq s(y)$, and consequently, $tx\in MS_1(A)$. The case of $MS_2(A)$ is proved similarly.

Now if $B=MS_1(A)\cup MS_2(A)$  
and $x\in B$ then $x\in MS_1(A)$ or $x\in MS_2(A)=-MS_1(A)$. If $x\in MS_1(A)$ and $0\leq t\leq 1$ then $tx\in MS_1(A)$ by what was proved above. 
If $-1\leq t\leq 0$ then $-tx\in MS_1(A)$, so $tx\in MS_2(A)$. In either case, $tx\in B$. The case $x\in MS_2(A)$ is similar.   
Therefore, $tx\in B$ whenever $x\in B$ and $|t|\leq 1$, and this shows that $B$ is balanced.
\end{bproof}

\begin{remark}
Regarding part (b) of the preceding proposition, it can be shown that if $A$ is absorbing then the assumption that the cone $\Vv_{sp}$ is generating is also a necessary condition for $MS_{1}(A)$ and $MS_{2}(A)$ to be absorbing. The proof of this fact can be based upon the properties of quasi-ideals in a mixed lattice space (see \cite[Theorem 4.12]{jj2}). However, we will not need this result in the present paper. 
\end{remark}

The next result provides additional conditions for our main characterization of mixed lattice topologies. As with the two types of mixed-full sets discussed earlier (see Proposition \ref{mixed-full-char}), we usually only need to work with one of the sets $MS_1(A)$ or $MS_2(A)$.

\begin{thm}\label{topology_charact}
Let $(\Vv,\tau)$ be a mixed lattice space with a generating cone $\Vv_{sp}$ and a vector topology $\tau$. 
Then the following conditions are equivalent. 
\begin{enumerate}[(a)]
\item
$\tau$ is locally mixed-full and the mixed lattice operations are continuous at zero.
\item
For every neighborhood $U$ of zero there exists a neighborhood $V$ of zero such that $V\subseteq U$ and $MS_1(V)\subseteq U$.
\item
For every neighborhood $U$ of zero there exists a neighborhood $V$ of zero such that $V\subseteq U$ and $MS_2(V)\subseteq U$. 
\end{enumerate}
\end{thm}

\begin{bproof}
\emph{(a)}$\implies$\emph{(b)} \, 
Let $U$ be a neighborhood of zero and let $W$ be a balanced mixed-full neighborhood of zero  
such that $W+W\subseteq U$. 
By the continuity of $x\mapsto \llow{x}$ and 
$x\mapsto s(x)$ at zero (see Remark \ref{huom1}) we can find a neighborhood $V$ of zero such that $x\in V$ implies $\llow{x}\in W$ and $\pm s(x)\in W$. 
Now if $-s(x)\sleq y\leq s(x)$ 
then $y\in W$ since $W$ is mixed-full.  
Hence, $MS_1(V)\subseteq W\subseteq U$. In particular, $-s(x)\sleq \rupp{x}\leq s(x)$ implies $\rupp{x}\in W$, 
and so $x=\rupp{x}-\llow{x}\in W+W\subseteq U$. Hence $V\subseteq U$, as required.

\emph{(b)}$\implies$\emph{(c)} \,
If $U$ is a neighborhood of zero and $W$ is a balanced neighborhood of zero such that $W\subseteq U$, then we can choose a neighborhood $V$ of zero such that $MS_1(V)\subseteq W$. Since $W$ is balanced we have $-MS_1(V)=MS_2(V)\subseteq W\subseteq U$. 
The converse implication \emph{(c)}$\implies$\emph{(b)} is similar.

\emph{(b)}$\implies$\emph{(a)} \, 
Let $U$ be a neighborhood of zero and let $V$ be a neighborhood of zero such that $MS_1(V)\subseteq U$. 
If $0\sleq y \leq x$ with $x\in MS_1(V)$ then 
$-s(z)\sleq x \leq s(z)$ for some $z \in V$. 
Hence $-s(z) \sleq y  \leq s(z)$, so $y\in MS_1(V)$. This shows that $MS_1(V)$ is a mixed-full neighborhood of zero, so $\tau$ is locally mixed-full by Proposition \ref{mixed-full-char}. 
Now, if $x\in V$ then $-s(x)\sleq \rupp{x} \leq s(x)$, and so $\rupp{x}\in MS_1(V)\subseteq U$. This shows that the map $x\mapsto \rupp{x}$ is continuous at zero, and hence the equivalence of $(a)$ and $(b)$ is proved. 
\end{bproof}


The notions of Riesz subspaces and ideals are in a central role in the theory of Riesz spaces. In the theory of mixed lattice spaces, the two partial orderings give rise to different types of ideals (see \cite{jj1,jj2}). 
A subspace $\Ss$ of a mixed lattice space $\Vv$ is a \emph{mixed lattice subspace} of $\Vv$ if the elements $x\uenv y$ and $x\lenv y$ are in $\Ss$ whenever $x,y\in\Ss$. 
A mixed-full mixed lattice subspace of $\Vv$ is called a \emph{quasi-ideal} of $\Vv$, and 
a $(\sleq)$-full mixed lattice subspace $\Aa$ of $\Vv$ is called a \emph{specific ideal} of $\Vv$. 
If $y,z\in \Vv$ then the sets $\{x\in \Vv: z\sleq x\leq y\}$ and $\{x\in \Vv: z\leq x\sleq y\}$ are called \emph{mixed-order intervals}.

We also recall that an ordered vector space $\Vv$ is called \emph{$(\leq)$-Archimedean} if $nx\leq y$ for all $n\in \N$ implies that $x\leq 0$.  
We shall also say that a sequence $\{x_n\}$ is \emph{$(\leq)$-increasing} if $x_n\leq x_m$ whenever $n\leq m$. A \emph{$(\sleq)$-increasing} sequence is defined similarly. 
The following theorem gives some additional properties of topological mixed lattice spaces. Similar results hold in topological Riesz spaces (see \cite{locsol}).

\begin{thm}\label{tms_prop}
Let $(V,\tau)$ be a topological mixed lattice vector space. Then the following statements hold.
\begin{enumerate}[(a)]
\item 
The positive cones $V_p=\{x\in V:x\geq 0\}$ and $V_{sp}=\{x\in V:x\sgeq 0\}$ are closed if and only if $\tau$ is a Hausdorff topology.
\item 
If $\tau$ is Hausdorff then $V$ is $(\leq)$-Archimedean. 
\item 
If $\tau$ is Hausdorff and $\{x_n\}$ is a $(\leq)$-increasing (or $(\sleq)$-increasing) sequence such that $x_n \stackrel{\tau}{\longrightarrow} x$ then $\sup \{x_n\}=x$ (or $\spsup \{x_n\}=x$, respectively).
\item 
If $\tau$ is locally mixed-full  and $A$ is a bounded set then $MF_1(A)$ and $MF_2(A)$ are bounded. 
\item 
If $\tau$ is locally mixed-full then 
every mixed-order interval is bounded. In particular, every ($\sleq$)-order interval is bounded. 
\item 
If $\Ss$ is a mixed lattice subspace of $\Vv$ then the closure $\overline{\Ss}$ is also a mixed lattice subspace. 
\item
The closure of a quasi-ideal is a quasi-ideal, and the closure of a specific ideal is a specific ideal.
\end{enumerate}
\end{thm}

\begin{bproof}
\begin{enumerate}[(a)]
\item
If $\tau$ is Hausdorff then the set $\{0\}$ is closed, and so the cones $V_p=\{x\in V:\rlow{x}= 0\}$ and $V_{sp}=\{x\in V:\llow{x}= 0\}$ are inverse images of a closed set under continuous mappings $x\mapsto \rlow{x}$ and $x\mapsto \llow{x}$, respectively, and are therefore closed. Conversely, if $V_p$ and $V_{sp}$ are closed, then $V_p\cap (-V_{sp})=\{0\}$ is closed and this implies that $\tau$ is Hausdorff.
\item
If $nx\leq y$ for all $n\in \N$ then $n^{-1}y-x\geq 0$ for all $n\in \N$, and so $-x=\lim_{n\to\infty} n^{-1}y-x \geq 0$ since $V_p$ is closed by (a). Hence $x\leq 0$ and $V$ is $(\leq)$-Archimedean.
\item
Let $\{x_n\}$ be a $(\leq)$-increasing sequence such that $x_n \stackrel{\tau}{\longrightarrow} x$. Then $x_n\leq x_{n+m}$ for all $m\in \N$, so $0\leq x_{n+m}-x_n$ and since by (a) $\Vv_p$ is closed, we have $\lim_{m\to\infty} (x_{n+m}-x_n)=x-x_n \geq 0$. Hence $x_n\leq x$ for all $n$. Suppose that $y$ is another upper bound of $\{x_n\}$. Then $y-x_n\geq 0$ for all $n$, and  $\lim_{n\to\infty} (y-x_n)=y-x \geq 0$, or $x\leq y$. This shows that $x=\sup \{x_n\}$. The case of a $(\sleq)$-increasing sequence is treated similarly.
\item
Let $A$ be a bounded set. If $U$ is any neighborhood of zero then there is an abrorbing mixed-full neighborhood $V$ of zero such that $V\subseteq U$. Since $A$ is bounded, there is some $\alpha \geq 0$ such that $A\subseteq \alpha V$. Then, since $\alpha V$ is mixed-full, we have $MF_1(A)\subseteq \alpha V\subseteq \alpha U$. This shows that $MF_1(A)$ is bounded. Similarly, $MF_2(A)$ is bounded. 
\item
Let $S=\{z:x\sleq z\leq y\}$ and let $V$ be a mixed-full and absorbing neighborhood of zero. Then there exist $\alpha >0$ such that $x\in \alpha V$ and $y\in \alpha V$. Then $S\subseteq \alpha V$, and so $S$ is bounded. A $(\sleq)$-order interval $\{z:x\sleq z\sleq y\}$ is contained in $S$, and is therefore bounded.
\item
Let $\Ss$ be a mixed lattice subspace of $\Vv$. Evidently, $\overline{\Ss}$ is a subspace. If $x\in \overline{\Ss}$ then there exists a net $\{x_\alpha\}\subset \Ss$ such that $x_{\alpha}\stackrel{\tau}{\longrightarrow} x$. Since $\Ss$ is a mixed lattice subspace, we have $\rupp{(x_\alpha)}\in \Ss$ for all $\alpha$, and $\rupp{(x_\alpha)}\stackrel{\tau}{\longrightarrow} \rupp{x}$ by the continuity of the map $x\mapsto \rupp{x}$. Hence 
$\rupp{x}\in \overline{\Ss}$ and so $\overline{\Ss}$ is a mixed lattice subspace. 
\item
Let $\Aa$ be a quasi-ideal. If $x\in \overline{\Aa}$ then there exists a net $\{x_{\alpha}\}\subset \Aa$ such that $x_{\alpha}\stackrel{\tau}{\longrightarrow} x$. If we  define another net $\{v_{\alpha}\}$ by ${v_{\alpha}}=x- {x_{\alpha}}$, then 
$v_{\alpha}\stackrel{\tau}{\longrightarrow} 0$. Now let $0\sleq y\leq x$. 
Since $x-v_{\alpha}\in \Aa$ it follows that $\rupp{(x-v_{\alpha})}\in \Aa$, and so $y-v_{\alpha}\leq x-v_{\alpha}$ for all $\alpha$. By \eqref{p4} this implies that $0\sleq \rupp{(y-v_{\alpha})}\leq \rupp{(x-v_{\alpha})}$, and consequently, $\rupp{(y-v_{\alpha})}\in \Aa$ for all $\alpha$. But $0\sleq y$, and by the continuity of the map $x\mapsto \rupp{x}$ we have $\rupp{(y-v_{\alpha})} \stackrel{\tau}{\longrightarrow} \rupp{y}=y$. Hence, $y\in \overline{\Aa}$ proving that $\overline{\Aa}$ is mixed-full. It now follows from part (f) that $\Aa$ is a quasi-ideal. The statement concerning specific ideals is proved in a similar way (just replace the assumption $0\sleq y\leq x$ by $0\sleq y\sleq x$ and apply \eqref{p6} in place of \eqref{p4}). 
\end{enumerate}
\end{bproof}


\section{Locally convex topologies on mixed lattice spaces}
\label{sec:sec35}

In a locally convex space the topology is determined by a family of seminorms, and in the case of a locally convex mixed lattice space the topology is compatible if the seminorms satisfy an additional condition.

\begin{defn}\label{seminorm}
Let $p$ be a seminorm in a mixed lattice vector space $\Vv$.
\begin{enumerate}[(i)]
\item
If $0\sleq x\leq y$  
implies $p(x)\leq p(y)$ then $p$ is called a \emph{mixed-monotone seminorm}.
\item
If $s(x)\leq s(y)$ implies $p(x)\leq p(y)$ then $p$ is called a \emph{mixed lattice seminorm}.  
\end{enumerate}
If the topology of $\Vv$ is given by a mixed lattice norm then $\Vv$ is called a \emph{normed mixed lattice space}.
\end{defn}

Locally convex mixed-full topologies have the following characterization.

\begin{thm}\label{loc_convex_mfull}
A locally convex topology $\tau$ on a topological mixed lattice vector space is locally mixed-full if and only if $\tau$ is generated by a family of mixed-monotone seminorms.
\end{thm}

\begin{bproof}
If $\tau$ is locally convex and mixed-full then by Proposition \ref{nbhoods} there exists a neighborhood base at zero $\{V_{\alpha}:\alpha\in I\}$ consisting of convex, absorbing, balanced, mixed-full sets. If $V_{\alpha}$ is any such neighborhood then the Minkowski functional $p_{\alpha}$ associated to $V_{\alpha}$ is a mixed-monotone seminorm. Indeed, let $p_{\alpha}(x)=\inf \{t>0:x\in tV_{\alpha}\}$ for each $x\in \Vv$. If $0\sleq y\leq x$   
then for every $\alpha\in I$ we have $x\in p_{\alpha} (x) V_{\alpha}$. Since $V_{\alpha}$ is mixed-full, so is $p_{\alpha} (x) V_{\alpha}$ and it follows that $y\in p_{\alpha} (x) V_{\alpha}$. By the definition of $p_{\alpha}$, this implies that $p_{\alpha} (y) \leq p_{\alpha} (x)$. Hence, 
$\{p_{\alpha}:\alpha\in I\}$ is a family of mixed-full seminorms generating $\tau$. 

Conversely, let $\{p_{\alpha}:\alpha\in I\}$ be a family of mixed-full seminorms generating $\tau$. For every finite subset $S$ of $I$ and $\varepsilon >0$ we put 
$V(S,\eps)=\{x:p_{\alpha}(x) < \varepsilon \textrm{ for all }\alpha\in S\}$. Then $0\sleq x\leq y$ and $y\in V(S,\eps)$ implies $x\in V(S,\eps)$, and so the sets $V(S,\eps)$ form a base of mixed-full neighborhoods of zero.
\end{bproof}

It follows from Definition \ref{seminorm} that the topology generated by a family of mixed lattice seminorms is locally symmetric-solid, and so the mixed lattice operations are continuous, by Theorem \ref{locsymmsol2}. The proof of the next result is similar to the latter part of the proof of Theorem \ref{loc_convex_mfull}.  

\begin{prop}\label{mlseminorm_topol}
If a vector topology $\tau$ on a mixed lattice vector space is generated by a family of mixed lattice seminorms then $\tau$ is locally symmetric-solid.
\end{prop}

We also have the following relationship between mixed lattice seminorms and mixed-monotone seminorms.

\begin{prop}\label{mlseminorm_char}
Let $p$ be a seminorm on $\Vv$. 
The following conditions are equivalent.
\begin{enumerate}[(a)]
\item
$p$ is a mixed lattice seminorm.
\item
$p$ is a mixed-monotone seminorm and $p(s(x))=p(x)$ for all $x\in \Vv$.
\end{enumerate}
\end{prop}

\begin{bproof}
If $(a)$ holds then, since $s(x)=s(s(x))$ for all $x\in \Vv$, it follows that $p(s(x))=p(x)$. If $0\sleq x\leq y$ then $0\sleq s(x)=\rupp{x}\leq \rupp{y}=s(y)$, and by assumption this implies that $p(x)\leq p(y)$.   
Conversely, if $(b)$ holds and $0\sleq s(x)\leq s(y)$ then $p(s(x))\leq p(s(y))$. By assumption, this is equivalent to $p(x)\leq p(y)$, and so $(a)$ follows.
\end{bproof}

As an immediate consequence of the preceding result, we notice  
that every mixed-monotone seminorm $p$ gives rise to a mixed-lattice seminorm $q$, defined by $q(x)=p(s(x))$ for all $x$, which we call a \emph{mixed lattice seminorm associated with} $p$.

\begin{prop}\label{ass_seminorm}
If $p$ is any mixed-monotone seminorm on $\Vv$ then $q(x)=p(s(x))$ is a mixed lattice seminorm. 
Moreover, if $p$ is a mixed-monotone norm, then $q$ is a mixed lattice norm.
\end{prop}

\begin{bproof}
If $0\sleq x\leq y$ then $0\sleq s(x)=\rupp{x}\leq \rupp{y}=s(y)$. Since $p$ is a mixed-monotone seminorm, it follows that $q(x)=p(s(x))\leq p(s(y))=q(y)$, so $q$ is also mixed-monotone. Moreover, $q(s(x))=p(s(s(x)))=p(s(x))=q(x)$, and so $q$ is a mixed lattice seminorm, by Proposition \ref{mlseminorm_char}. 
If $p$ is a norm then so is $q$, because then $q(x)=p(s(x))=0$ iff $s(x)=0$ iff $x=0$.
\end{bproof}

With the concepts and results introduced above we can now complete our characterization of mixed lattice topologies. It was shown in Theorem \ref{locsymmsol2} that a locally symmetric-solid topology is locally mixed-full. If $\Vv$ is locally convex then we have the following converse.

\begin{thm}\label{loc_convex_char1}
Let $(\Vv,\tau)$ be a locally convex topological mixed lattice space such that the cone $\Vv_{sp}$ is generating. 
If $\tau$ is locally mixed-full then $\tau$ is locally symmetric-solid.
\end{thm}

\begin{bproof}
If $\tau$ is a locally convex and locally mixed-full then $\tau$ is determined by a directed family $\{p_{\alpha}:\alpha\in I\}$ of mixed-monotone seminorms, by Theorem \ref{loc_convex_mfull}. 
For each $p_{\alpha}$, let $q_{\alpha}$ be the associated mixed lattice seminorm, that is, $q_{\alpha}(x)=p_{\alpha}(s(x))$ for all $x\in \Vv$, and let us denote the topology determined by the family $\{q_{\alpha}:\alpha\in I\}$ by $\widetilde{\tau}$. By Proposition \ref{mlseminorm_topol} $\widetilde{\tau}$ is locally symmetric solid. We will show that $\tau=\widetilde{\tau}$.

For every $\varepsilon >0$ and $\alpha\in I$ the set $W=\{x:p_{\alpha}(x) < \frac{\eps}{2}\}$ is a mixed-full $\tau$-neighborhood of zero, and $\widetilde{W}=\{x:q_{\alpha}(x) < \eps \}$ is a symmetric-solid $\widetilde{\tau}$-neighborhood of zero. 
%
By Theorem \ref{topology_charact} there exists a 
$\tau$-neighborhood $W_0\subseteq W$ of zero such that $MS_1(W_0)\subseteq W$. Thus, for every $x\in MS_1(W_0)$ there exists $y\in W_0$ such that $-s(y)\sleq x \leq s(y)$ and $s(y)\in MS_1(W_0)$. 
This inequality implies that $0\sleq \rupp{x} \leq s(y)$ and $0\sleq \rlow{x} \leq s(y)$, and summing these inequalities gives $s(x)\leq 2s(y)$. This implies that $q_{\alpha}(x)\leq 2q_{\alpha}(y)=2p_{\alpha}(s(y)) < \eps$. Hence $x\in \widetilde{W}$, and so $MS_1(W_0) \subseteq \widetilde{W}$, proving that $\widetilde{\tau}\subseteq \tau$.

For the converse, choose the neighborhoods $V=\{x:p_{\alpha}(x) < \eps\}$  and $\widetilde{V}=\{x:q_{\alpha}(x) < \frac{\eps}{3}\}$. Then, in particular, $\widetilde{V}$ is mixed-full since $q_{\alpha}$ is a mixed-monotone seminorm, by Proposition \ref{mlseminorm_char}. Thus, there exists a 
$\widetilde{\tau}$-neighborhood $\widetilde{W}_0\subseteq \widetilde{W}$ 
of zero such that $MS_1(\widetilde{W}_0)\subseteq \widetilde{W}$. Thus, for every $x\in MS_1(\widetilde{W}_0)$ there exists $y\in \widetilde{W}_0$ such that $-s(y)\sleq x \leq s(y)$ and $s(y)\in MS_1(\widetilde{W}_0)$. Then $0\sleq x+s(y) \leq 2s(y)$, and since $p_{\alpha}$ is a mixed-monotone seminorm, we get $p_{\alpha}(x+s(y))\leq 2p_{\alpha}(s(y)$. The triangle inequality then gives $p_{\alpha}(x)-p_{\alpha}(s(y))\leq 2p_{\alpha}(s(y))$, or $p_{\alpha}(x)\leq 3p_{\alpha}(s(y))=3q_{\alpha}(x) < \eps$. Therefore, $x\in W$, and so $MS_1(\widetilde{W}_0)\subseteq W$. This shows that $\tau\subseteq \widetilde{\tau}$, and hence $\tau = \widetilde{\tau}$.
\end{bproof}

Combining the results of the preceding theorem with Theorem \ref{locsymmsol2} 
gives a complete characterization of mixed lattice topologies in the locally convex case, assuming the cone $\Vv_{sp}$ is generating. 

\begin{cor}\label{complete}
Let $\Vv$ be a mixed lattice space such that the cone $\Vv_{sp}$ is generating. If $\tau$ is a locally convex topology on $\Vv$ then the following conditions are equivalent.
\begin{enumerate}[(a)]
\item
The mixed lattice operations are uniformly continuous.
\item
$\tau$ is locally symmetric-solid.
\item
$\tau$ is locally mixed-full and the mixed lattice operations are continuous at zero.
\item
$\tau$ is locally ($\leq$)-full and the mixed lattice operations are continuous at zero.
\end{enumerate} 
\end{cor}

Next we investigate the relationships between mixed lattice topologies and locally solid Riesz space topologies. If $\Vv$ is a mixed lattice vector space which is a lattice with respect to the specific order $\sleq$ then $\Vv$ can be equipped with a locally solid (with respect to $\sleq$) Riesz space topology, and such topology will be called a \emph{locally $(\sleq)$-solid} Riesz space topology. Similarly, if $\Vv$ is a lattice with respect to $\leq$ then we use the term \emph{locally $(\leq)$-solid} Riesz space topology. As usual, if $(\Vv,\leq)$ is a Riesz space then the lattice-theoretic positive part  
of an element $x$ is denoted by $x^+$.  
Similarly, if $(\Vv,\sleq)$ is a Riesz space then the positive part  
of $x$ with respect to $\sleq$ is denoted by $\spec (x^+)$,  
and the least upper bound of elements $x$ and $y$ with respect to $\sleq$ is denoted by $\spsup \{x,y\}$.

\begin{thm}\label{splattice}
Let $(\Vv,\leq,\sleq)$ be a mixed lattice vector space which is a lattice with respect to specific order $\sleq$. If $\tau$ is a locally $(\sleq)$-solid Riesz space topology on $\Vv$ then the mixed lattice operations are continuous.
\end{thm}

\begin{bproof}
In \cite[Proposition 3.16]{jj1} it was shown that $\spec (x^+)=\spsup \{\lupp{x},\rupp{x}\}$, and so we have $0\sleq \rupp{x} \sleq \spec (x^+)$. Now if $U$ is a neighborhood of zero and $\tau$ is locally $(\sleq)$-solid, we can choose a $(\sleq)$-solid neighborhood $W$ of zero such that $W\subset U$, and by the continuity of the mapping $x\mapsto \spec (x^+)$ there exists a neighborhood $V$ of zero such that $x\in V$ implies $\spec (x^+)\in W$. Hence, by the
above inequality $x\in V$ implies $\rupp{x}\in W\subset U$, and hence the map $x\mapsto\rupp{x}$ is continuous at zero. By Proposition \ref{unicont}, all the mixed lattice operations are then continuous.
\end{bproof}

\begin{thm}\label{initlattice}
Let $(\Vv,\leq,\sleq)$ be a  
mixed lattice vector space which is a lattice with respect to initial order $\leq$. If $\tau$ is a locally symmetric-solid topology on $\Vv$ then $\tau$ is a locally $(\leq)$-solid Riesz space topology.
\end{thm}

\begin{bproof}
By Theorem \ref{locsymmsol2} $\tau$ is locally $(\leq)$-full and the map $x\mapsto\lupp{x}$ continuous at zero. By the Roberts--Namioka Theorem \ref{rob_nam} it is sufficient to show that the map $x\mapsto {x}^+$ is continuous at zero. To this end, we just need to note that $\lupp{x}\geq x$ and $\lupp{x}\geq 0$, and hence $\lupp{x}\geq x^+$. Now if $U$ is a neighborhood of zero, we can choose a $(\leq)$-full neighborhood $V$ of zero such that $V\subset U$, and by the continuity of $x\mapsto\lupp{x}$ there is a neighborhood $W$ of zero such that $x\in W$ implies $\lupp{x}\in V$. But $0\leq x^+ \leq \lupp{x}$ and $V$ is full, so $x^+ \in V\subset U$ and this completes the proof.
\end{bproof}

As for the lattice properties of topological mixed lattice spaces, our results have the consequence that finite-dimensional normed mixed lattice spaces with a generating cone $\Vv_{sp}$ are necessarily lattices with respect to specific order $\sleq$. To prove this, we need the following theorem which was proved by Arsove and Leutwiler for mixed lattice semigroups (\cite[Theorem 9.7]{ars}), and essentially the same proof carries over to mixed lattice vector spaces (cf. \cite[Theorem 3.9]{eri}). 

\begin{thm}[{\cite[Theorem 3.9]{eri}}]\label{splattice}
Let $\Vv$ is a mixed lattice space such that the cone $\Vv_{sp}$ is generating and $\spsup \{x_n\}$ exists for every $(\sleq)$-increasing and $(\sleq)$-bounded sequence $\{x_n\}$. Then $\Vv$ is a lattice with respect to $\sleq$.
\end{thm}

The following theorem is now rather immediate.

\begin{thm}\label{finitedimlattice}
Let $\Vv$ be a finite-dimensional normed mixed lattice space. If the cone $\Vv_{sp}$ is generating then $\Vv$ is a lattice with respect to $\sleq$.
\end{thm}

\begin{bproof}
Let $\{x_n\}$ be a $(\sleq)$-increasing sequence such that $x_n \sleq u$ for all $n$ and some $u\in \Vv$. Since $\{x_n\}$ is contained in the $(\sleq)$-order interval $\{z:x_1\sleq z \sleq u\}$, it follows by Theorem \ref{tms_prop}(e) that $\{x_n\}$ is norm bounded. As $\Vv$ is finite-dimensional and $\{x_n\}$ is $(\sleq)$-increasing, it follows from the Bolzano--Weierstrass theorem that $\{x_n\}$ converges to a limit $x$, and by Theorem \ref{tms_prop}(c) we have $x=\spsup \{x_n\}$. Now Theorem \ref{splattice} implies that $\Vv$ is a lattice with respect to $\sleq$.   
\end{bproof}


Locally solid Riesz spaces are of course a special case of topological mixed lattice spaces. 
Let us consider the following example to illustrate our results.

\begin{example}\label{bv}
Let $BV([a,b])$ be the space of functions of bounded variation on an interval $[a,b]$.  
We define initial order in $BV([a,b])$ by
$$
f\leq g \qquad \iff \qquad f(x)\leq g(x) \quad \textrm{for all } x\in[a,b]
$$
and specific order by
$$
f\sleq g \; \iff \; f(x)\leq g(x) \; \textrm{for all } x\in[a,b] \; \textrm{and} \; g-f \; 
\textrm{ is increasing on } [a,b].
$$
It has been shown in \cite{eri} that $\Vv=(BV([a,b]),\leq,\sleq)$ is a mixed lattice vector space 
where the mixed lower and upper envelopes are given by 
$$
(f\lenv g)(u)=\inf \, \{f(u)- (f(x)-g(x))^+ \, : \, x\in[a,u] \}
$$
and
$$
(f\uenv g)(u)=\sup \, \{f(u) + (g(x)-f(x))^+ \, : \, x\in[a,u] \},
$$
where $c^+ = \max \{0,c\}$ is the positive part of the real number $c$. 
It is well-known that $\Vv$ is a lattice with respect to both partial orderings $\leq$ and $\sleq$. In particular, the positive cone $\Vv_{sp}$ is generating. 

First we will show that if $\Vv$ is equipped with the sup-norm $\norm{f}_\infty =\sup \{|f(x)|:x\in [a,b]\}$ then $\Vv$ is a topological mixed lattice space. By Proposition \ref{mlseminorm_topol} and Corollary \ref{complete} it is sufficient to show that $s(f)\leq s(g)$ implies $\norm{f}_\infty\leq \norm{g}_\infty$.  
Let $g\in V$. Let us first find $s(g)$ using the above formulae for the mixed envelopes. We have
$$
\rupp{g}(u)=(0\uenv g)(u)=\sup_{x\in[a,u]} \, \{(g(x))^+ \},
$$
and
$$
\rlow{g}(u)=(0\uenv (-g))(u)=\sup_{x\in[a,u]} \, \{(-g(x))^+ \} = \sup_{x\in[a,u]} \, \{(g(x))^- \}.
$$
Since $s(g)=\rupp{g}+\rlow{g}$, we obtain
$$
s(g)(u)=
\sup_{x\in[a,u]} \, \{(g(x))^+ + (g(x))^-\}=\sup_{x\in[a,u]} \, \{|g(x)| \}.
$$
Hence we have
$$
\norm{g}_{\infty}=\sup \, \{|g(x)| \, : \, x\in[a,b] \}=s(g)(b).
$$ 
This shows that if $s(f)\leq s(g)$ then $\norm{f}_{\infty}=s(f)(b)\leq s(g)(b)=\norm{g}_{\infty}$, and so $\norm{\cdot }_{\infty}$ is a mixed lattice norm on $\Vv$ which generates a mixed lattice topology (by Proposition \ref{mlseminorm_topol}). In this case, this topology is the same as the locally solid Riesz space topology on $(\Vv,\leq)$ induced by the $\sup$-norm.

On the other hand, if $\Vv$ is equipped with the total variation norm $\norm{f}_{\scriptscriptstyle BV}=|f(a)|+V_a^b(f)$, (where $V_a^b(f)$ is the total variation of $f$ on $[a,b]$) then the mixed lattice operations are continuous. This follows from Theorem \ref{splattice}, since $\Vv$ is a vector lattice with respect to $\sleq$ and $\norm{\cdot }_{\scriptscriptstyle BV}$ is a $(\sleq)$-lattice norm. Moreover, if $0\sleq f \leq g$ then $f$ is positive and non-decreasing, and $f(x)\leq g(x)$ for all $x\in [a,b]$, so we have $\norm{f}_{\scriptscriptstyle BV}=f(b)\leq g(b)\leq \norm{g}_{\scriptscriptstyle BV}$. This shows that the topology generated by the total variation norm is locally mixed-full, and so $(\Vv,\norm{\cdot }_{\scriptscriptstyle BV})$ is a locally symmetric-solid mixed lattice space also with respect to the $BV$-norm topology, by Corollary \ref{complete}. Consequently, $(\Vv,\norm{\cdot }_{\scriptscriptstyle BV})$ is also a locally $(\leq)$-solid Riesz space, by Theorem \ref{initlattice}. 
Note, however, that the norms $\norm{\cdot }_\infty$ and $\norm{\cdot }_{\scriptscriptstyle BV}$ are not equivalent. 
\end{example}

\section{On asymmetric norms and cone norms}
\label{sec:sec4}

The theory of asymmetric normed spaces is relatively recent. A detailed account of the theory is given in the monograph \cite{cob}. Although we shall not explore these ideas in detail in the present paper, the discussion in the preceding sections 
suggests that topologies determined by asymmetric seminorms arise rather naturally in mixed lattice spaces. For example, we have seen that there are two types of mixed-full sets that behave in an asymmetric fashion, that is, if $A\subseteq \Vv$ is a balanced set then $MF_1(A)=-MF_2(A)$. On the other hand, if a mixed lattice norm on $\Vv$ is given, then the properties of the upper and lower parts give rise to asymmetric norms on $\Vv$, as the following discussion shows. 
The definitions and notations that follow are adopted from \cite{cob}.

Let $X$ be a real vector space. A real-valued mapping $p:X\to [0,\infty)$ is called an \emph{asymmetric norm} on $X$ if for all $x,y\in X$ the following conditions hold.

$$ 
\begin{array}{ll}
\text{(A1)} & \quad p(x)=0 \quad \text{and} \quad p(-x)=0 \quad \text{implies} \quad x=0 \\
\text{(A2)} & \quad p(\alpha x)=\alpha p(x) \quad \text{for all} \quad \alpha \geq 0 \\
\text{(A3)} & \quad p(x+y)\leq p(x)+p(y) \\[2ex]
\end{array}
$$

If $p$ satisfies conditions (A2) and (A3) then $p$ is called an \emph{asymmetric seminorm}. \\

Given an asymmetric norm $p$, its \emph{conjugate asymmetric norm} $\overline{p}$ is defined as $\overline{p}(x)=p(-x)$. If we put $p^s(x)=p(x)+\overline{p}(x)$ 
then $p^s$ is a norm associated with the asymmetric norm $p$.

To see how these concepts are related to the results of the preceding sections,  
we note that if $\Vv$ is a mixed lattice space and $U$ is a convex, absorbing and balanced neighborhood of zero then by Proposition \ref{nbhoods} the sets $MF_1(U)$ and $MF_2(U)$ are also convex and absorbing neighborhoods of zero (but not balanced, since $x\in MF_1(U)$ implies $tx\in MF_1(U)$ only for $0\leq t\leq 1$, and similarly for $MF_2(U)$).  
Hence, the Minkowski functionals $p_1$ of the set $MF_1(U)$ and $p_2$ of the set $MF_2(U)$ are asymmetric seminorms. In fact, they are conjugate asymmetric seminorms since $MF_1(U)=-MF_2(U)$. Furthermore, by Proposition \ref{nbhoods}, the sum $MF_1(U)+MF_2(U)$ is a balanced, convex, absorbing and mixed-full neighborhood of zero, whose Minkowski functional is a mixed-monotone seminorm by Theorem \ref{loc_convex_mfull}.

In normed Riesz spaces, asymmetric norms 
are related to lattice norms in a quite natural way through the positive and negative parts of an element.  Similarly, in mixed lattice spaces the upper and lower parts give rise to asymmetric norms. 
If $\rho$ is any mixed lattice norm on $V$ then we can define $p_1(x)=\rho(\rupp{x})$ and $p_2(x)=\rho(\rlow{x})$. Then $p_1$ and $p_2$ are conjugate asymmetric seminorms on $V$. It follows immediately from the properties of $\rupp{x}$ and $\rlow{x}$ (Theorem \ref{absval}) that $p_1$ and $p_2$ satisfy the conditions in the definition of an asymmetric norm. 



Asymmetric cone norms are further generalizations of asymmetric norms, and they have found applications in different areas of mathematics. 
In recent years, they have been used for studying generalizations of asymmetric normed spaces and some related results in analysis \cite{ilk, ilk1}. The basic idea in such generalizations is that the usual asymmetric norm is replaced by a vector-valued mapping that has properties similar to those of an asymmetric norm, but its range is a positive cone of a partially ordered vector space.  
Asymmetric cone norms are also closely related to certain problems in optimization theory \cite{nemeth2020}.

\begin{defn}\label{asymm_cone_norm}
Let $C\subset X$ be a cone in a topological vector space $X$ and let $\leq$ be the associated order relation. A continuous mapping $Q:X\to C$ is an \emph{asymmetric cone norm} if 
\begin{enumerate}[(1)]
\item 
$Q(x)=x$ for all $x\in C$ and $Q(X)=C$
\item
$Q(tx)=tQx$ for all $t\in \R_+$ and $x\in X$
\item
$Q(x+y)\leq Qx+Qy$ for all $x,y\in X$
\item
If $Qx=0$ and $Q(-x)=0$ then $x=0$.
\end{enumerate}
Moreover, we say that $Q$ is a \emph{proper asymmetric cone norm} if $Q(I-Q)=0$. Here $I$ denotes the identity operator on $X$.
\end{defn}

In a general topological vector space, the existence of a mapping satisfying the conditions in the above definition is not at all clear. However, the properties of the mixed lattice operations given in Theorem \ref{absval} yield the following existence result.

\begin{thm}\label{thm2}
Let $\Vv=(\Vv,\leq,\sleq)$ be a topological mixed lattice space. Then the following hold.
\begin{enumerate}[(a)]
\item
The mapping $Q:\Vv\to \Vv_{p}$ given by $Q(x)=\lupp{x}$ is a proper asymmetric cone norm on $\Vv$. 
\item
The mapping $Q:\Vv \to \Vv_{sp}$ given by $Q(x)=\rupp{x}$ is a proper asymmetric cone norm on $\Vv$ which is increasing with respect to both partial orderings (that is, $x\sleq y$ implies $Qx \sleq Qy$ and $x\leq y$ implies $Qx \leq Qy$).
\end{enumerate} 
\end{thm}

\begin{bproof}
\emph{(a)}\, The properties of $\lupp{x}$ imply 
that the mapping $Q(x)=\lupp{x}$ has the properties listed in Definition \ref{asymm_cone_norm}. Property $(1)$ follows from Theorem \ref{absval}(g), and $(2)$ follows from \eqref{p8a}. Property $(3)$ is an immediate consequence of Theorem \ref{absval}(c), while $(4)$ follows from Theorem \ref{absval}(a) and (g). 
To check that $Q$ is proper, we note that 
$$
Q(I-Q)(x)=Q(x-\lupp{x})=Q(-\rlow{x})=\lupp{(-\rlow{x})} =\llow{(\rlow{x})} =0,
$$ 
by Theorem \ref{absval} (since $\rlow{x}\sgeq 0$). 

The proof of \emph{(b)} is similar. 
The property that $Q$ is increasing with respect to both partial orderings 
follows by \eqref{p4} and \eqref{p6}.  
\end{bproof}

Suppose that $C$ is a cone in a topological vector space $X$, and we ask if there exists an asymmetric cone norm $Q$ on $X$ 
satisfying the conditions in Definition \ref{asymm_cone_norm}. 
Theorem \ref{thm2} shows that a sufficient condition for the existence of such mapping is that there exists a mixed lattice order structure on $X$ such that $C$ is a positive cone for the partial order $\leq$.  

As an application of our results, we show that given any closed convex cone $C$ in $\R^n$ we can always turn $\R^n$ into a topological mixed lattice space $(\R^n,\leq,\sleq)$ in such way that $C$ will be the positive cone associated with the partial order $\leq$, and so by Theorem \ref{thm2} there always exists a continuous asymmetric cone norm associated with the given cone in $\R^n$. The latter part (that is, the existence of a proper asymmetric cone norm) has been proved recently in \cite[Theorem 2]{nemeth2020} using different methods. This result is of importance in problems related to optimization theory (see \cite{nemeth2020} and references therein).


\begin{thm}\label{thm1}
Let $C$ be a generating closed convex pointed cone in $\R^n$ 
and let $x\in \interior (C)$. If $\leq$ is the partial ordering induced by $C$ and $\sleq$ is another partial ordering induced by the cone $R_x=\{t x:t \geq 0\}$ then $(\R^n,\leq,\sleq)$ is a normed mixed lattice vector space with respect to the usual norm of $\R^n$. In particular, the mapping $Q:\R^n\to C$ given by $Q(x)=\lupp{x}$ is a continuous proper asymmetric cone norm. 
\end{thm}

In the proof of Theorem \ref{thm1} we use the well known fact that the interior $\interior (C)$ of a generating cone $C$ in $\R^n$ is non-empty. 
We also 
recall that a subset $B$ of a cone $C$ is called a \emph{base of} $C$ if for every $x\in C\setminus \{0\}$ there exists a unique number $t>0$ such that $tx\in B$. 
We will also need the following geometric result, which is a consequence of the existence of a hyperplane supporting $C$ at $0$ (c.f. \cite[Lemma 2]{nemeth2020}). 

\begin{lemma}
If $C\subset \R^n$ is a closed convex cone such that $\interior (C)\neq \emptyset$ then there exists a closed, bounded and convex base $B$ of $C$ and an interior point $x\in \interior (C)$ such that the line spanned by $x$ is orthogonal to $B$.
\end{lemma}

\begin{bproof}[Proof of Theorem \ref{thm1}]
Let $\leq$ be the partial ordering induced by the cone $C$. 
By the preceding lemma, there is a closed, bounded and convex base $B$ of $C$ and an interior point $x$ of $C$ such that $B$ is orthogonal to ray $R_x$ generated by $x$ (i.e. $R_x=\{tx:t\geq 0\}$). 
Clearly, $R_x$ is a cone and it induces a partial ordering $\sleq$ on $\R^n$ by $u\sleq v$ if $v-u\in R_x$, or equivalently, $v=u+tx$ for some $t\geq 0$. Evidently, $R_x$ is also a mixed lattice cone, as it is ''one dimensional'' and both partial orderings coincide on $R_x$. Now $V=(\R^n,\leq,\sleq)$ is a partially ordered vector space with two partial orderings such that $u\sleq v$ implies $u\leq v$. 
To show that $V$ is a mixed lattice space it is sufficient to show that the element 
$$
\lupp{y}=y\uenv 0 = \min \,\{\,w\in V: \; w\sgeq y \; \textrm{ and } \; w\geq 0 \,\}
$$ 
exists for all $y\in V$ (by Theorem \ref{mlg_ehto}). We note that in the present situation this is equivalent to showing that 
$$
\min \,(\{y+tx :t\geq 0 \} \cap C )= \min \{t\geq 0: y+tx\in C\}
$$
exists for all $y\in V$. 
If we denote the boundary of $C$ by $\partial C$ then the distance of $x$ to $\partial C$ is defined as $d(x,\partial C)=\inf\{\norm{x-y}:y\in \partial C\}$. 
Then the set $D$ bounded by $B$ and $\partial C$ and containing $0$ is convex, closed and bounded. 
Moreover, $tx\in B$ for some $t>0$, so we may thus assume that $x\in D$. 
%
Now $d (x,\partial C)>0$ since $x$ is an interior point of $C$, and since $D$ is closed, convex and bounded there is a point $v\in D\cap\partial C$ (not unique, in general) such that $d (x,v)=d (x,D\cap\partial C)=d (x,\partial C)$. The last equality follows from the orthogonality of $B$ and $R_x$. 
Hence, we have shown that for any $x\in \interior (C)$ there is some point $v\in \partial C$ such that $d(x,\partial C)=d(x,v)$.  
Now for any $t>0$ we have $d (tx,\partial C)=d (tx,tv)=td (x,v)$. 
Indeed, if there was some $u\in \partial C$ such that $\norm{tx-u}< \norm{tx-tv}=t\norm{x-v}$, then this would imply that $\norm{x-t^{-1}u}< \norm{x-v}$, a contradiction. 
Therefore, if $0<s <t$ then 
$d (s x,\partial C) 
< d (t x,\partial C)$,   
and so the distance of $tx$ to the boundary of $C$ increases 
as $t$ increases. But $\{tx:t\geq 0\}$ and $\{y+tx:t\geq 0\}$ are parallel half-lines, so $d (t x,y+tx)=\norm{y}$ is a constant. Thus, there is some $t_1$ such that $d (t x,\partial C)\geq d (t x,y+tx)=\norm{y}$ for all $t\geq t_1$, and so $y+tx$ lies inside $C$ for all $t\geq t_1$. Since $C$ is closed, the number $t_0=\min \{t\geq 0:y+tx\in C\}$ exists, and 
so the element $\lupp{y}=y+t_0 x$ exists for all $y\in \R^n$. This shows that $(\R^n,\leq,\sleq)$ is a mixed lattice vector space, by Theorem \ref{mlg_ehto}. Then by Proposition $\ref{thm2}$ the mapping $Q(x)=\lupp{x}$ is a proper asymmetric cone norm.

We still need to prove the continuity of $Q$. Let $\norm{\cdot }_2$ be the usual Euclidean norm on $\R^n$. Define a new norm by $\norm{z}_0=\norm{s(z)}_2$ for all $z\in \R^n$. We will show that $\norm{\cdot }_0$ is a mixed lattice norm. Assume that $s(x)\leq s(y)$. Then, since $s(x),s(y)\sgeq 0$, it follows that $s(x)$ and $s(y)$ are both contained in the ray $R_x$. Hence, $s(x)\leq s(y)$ implies that $s(y)=t s(x)$ for some $t\geq 1$, and so $\norm{s(x)}_2\leq t\norm{s(x)}_2=\norm{s(y)}_2$, or $\norm{x}_0 \leq \norm{y}_0$. This shows that $\norm{\cdot }_0$ is a mixed lattice norm, and hence the mapping $Q$ is continuous in the topology generated by the norm $\norm{\cdot }_0$. Since all norms on $\R^n$ are equivalent, it follows that the mapping $Q$ is continuous in the usual topology of $\R^n$.
\end{bproof}




We remark that the proof of Theorem \ref{thm1} is based on rather simple geometric arguments which do not necessarily work in infinite-dimensional spaces. For instance, a generating cone in an infinite dimensional space does not always possess interior points. 



\bibliographystyle{plain}

\end{document}